\documentclass[11pt]{article}
\usepackage{amsfonts,amsmath,
amscd,latexsym,color,epsfig,amsthm}
\usepackage{epsfig}

    \textwidth 
    12.5cm   
    \textheight 
    18.5cm

\date{}

\title{A colouring problem for the dodecahedral graph}

\author{
{\sl Endre Makai, Jr.,}\thanks{Alfr\'ed R\'enyi Institute of Mathematics,
Hungarian Academy of Sciences, H-1364 Budapest, P.O. Box 127, Hungary;
\texttt{makai.endre@renyi.mta.hu}; www.renyi.mta.hu/\~{}makai.}
{\sl Tibor Tarnai}{\thanks{Budapest University of Technology and Economics,
Department of Structural Mechanics, H-1521 Budapest, M\H{u}egyetem rkp.\ 3, 
Hungary; \texttt{tarnai@ep-mech.me.bme.hu}; www.me.bme.hu/tarnai-tibor.
Partially supported by Hungarian National Research, Development and Innovation
Office, NKFIH, Grant No. K119440.}}}

\begin{document}

\maketitle



\medskip


We investigate an elementary colouring problem for the dodecahedral graph.
We determine the number of all colourings satisfying a certain condition. We
give a simple combinatorial proof, and
also a geometrical proof. For the proofs we use the symmetry group of the
regular dodecahedron, and also the compounds of five tetrahedra, inscribed in
the dodecahedron. Our results put a result of W. W. Rouse Ball -- H. S. M.
Coxeter in the proper interpretation.

{\bf{CV-s.}} 

E. Makai, Jr. graduated in mathematics at L. E\"otv\"os University,
Budapest in 1970.
He is a professor emeritus of Alfr\'ed R\'enyi Institute of Mathematics,
Hungarian Academy of Sciences. His main interest is geometry.

T. Tarnai graduated in civil engineering 
at the Technical University of Budapest in 1966,
and in mathematics at L. E\"otv\"os University, Budapest in 1973.
He is a professor emeritus of Budapest University of Technology and
Economics. His main interests are engineering and mathematics, and in
mathematics connections between engineering methods and geometry.


\newpage


\section{Introduction}

As well known, the most popular soccer balls, i.e.,\ footballs, 
have the form of a
spherical mosaic, consisting of twelve regular spherical
pentagons and twenty regular spherical 
hexagons. If we retain all the vertices of all its faces, but replace the
regular spherical polygons by planar ones, then we obtain an 
{\it{Archimedean polyhedron}}. That is, it is convex,
its faces are regular, and its symmetry
--- i.e., congruence --- 
group acts transitively on its vertices. The group of symmetries of this
spherical mosaic coincides with the group of symmetries of this Archimedean
polyhedron, and also with the group of symmetries of the regular dodecahedron.
Since the regular dodecahedron is the polar reciprocal 
of the regular icosahedron, 
the symmetry group of the regular dodecahedron coincides with that of the
regular icosahedron. Therefore this group is called the
{\it{icosahedral group}}, denoted by $I_h$. Its
elements preserve the centre of the regular dodecahedron.
Then $I$ denotes the subgroup of $I_h$,
consisting of the orientation preserving elements of $I_h$. Then
$I_h$ has $120$, and $I$ has $60$ elements. Their descriptions will be given 
in {\bf{1.1}}.
Both the above spherical
mosaic, and the corresponding Archimedean polyhedron are denoted by $(5,6,6)$,
cf.\ [1]. The numbers $5,6,6$ refer to the fact that at any vertex there
meet one pentagon and two hexagons.

This spherical mosaic
was invented as a soccer ball by the former Danish football player Eigil
Nielsen and was introduced in 1962.
The usual colouring of the soccer balls is the following:
the pentagons are coloured black, and the hexagons white, patented by M. Doss
in 1964.

However, this mosaic was used for making balls much earlier. An 
archive home movie made at the end of the 1930s shows the appearance of such
a ball, which we call the W\"orthsee ball.
However the name of the maker of this ball and the year of its making are
unknown to us.
The faces of this ball were coloured as follows: 
all pentagons were green, while the
hexagons were of five different colours. In our notation in Fig.~1,
white~=~1, yellow~=~2, red~=~3, blue~=~4, black~=~5.
Seemingly the intention of the maker was that each pentagon should be
surrounded by hexagons of all five colours.
This was achieved for several, but not all pentagons, cf.\ Fig.~1(A).
For a more detailed history cf. T. Tarnai -- A. Lengyel \cite{TL}.




\smallskip
\vbox{\centerline{\epsfxsize=61mm 
\epsfbox{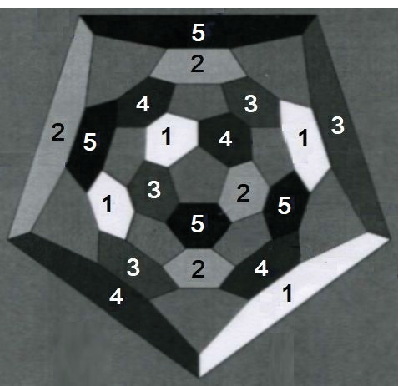}\hfill
{\epsfxsize=64mm \epsfbox{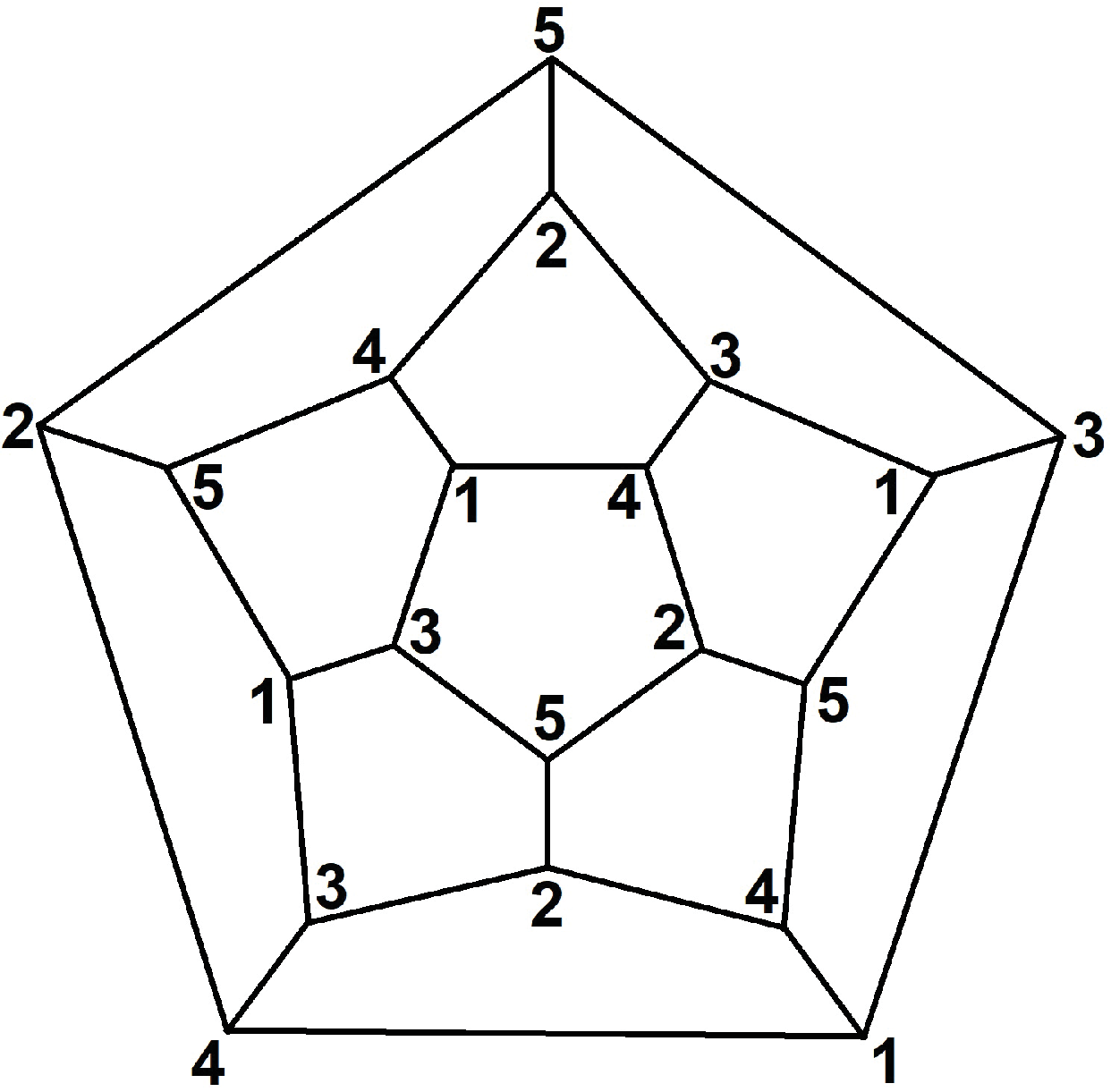}}}

\centerline{\hfil(A)\hfil\hfil\hfil\hfil(B)\hfil}

\smallskip
\noindent
{
Figure 1.
Colouring of the W\"orthsee ball. (A)~Schlegel diagram of $(5,6,6)$ with 
coloured faces.
(B)~Schlegel diagram of the corresponding regular dodecahedron with coloured
vertices.
\par}}

\smallskip


This raises the following mathematical problem.
Does there exist a colouring of the hexagonal faces of $(5,6,6)$ with five
colours, such that each pentagon has hexagonal neighbours of all five colours?
Clearly in this problem the size of the pentagons is not important:
we may shrink them simultaneously to their centres, while the hexagonal faces
become regular triangles.
Then our problem becomes the following.
Colour the twenty faces of a regular icosahedron with five colours, so that each
vertex is incident to faces of all five colours.

It will be more convenient to consider the dual problem.
Colour the twenty vertices of the regular dodecahedron so that each face
should have vertices of all five colours.

The actual colouring of the dual of the
W\"orthsee ball, with the pentagons of the W\"orthsee ball
being shrunken to their centres, is shown in Fig.~1(B).


This colouring problem was considered later in mathematics as well.

H. M. Cundy -- A. P. Rollett, \cite{CR}, 1985, write on pp.\ 82--83 the
following.
``The faces of the icosahedron can be coloured by five colours so that the
five faces at every vertex are coloured differently, but opposite faces cannot
then be coloured alike.''

W. W. Rouse Ball -- H. S. M. Coxeter, \cite{RBC}, 1987, write on p.\ 242 the
following.
``L. B. Tuckerman remarks that the faces of an icosahedron can be coloured with
five colours so that each face and its three neighbours have four different
colours. For five given colours
this can be done in four ways, consisting of two enantiomorphous pairs.
One pair can be derived from the other by making any odd permutation of the
five colours, e.g.,\ by interchanging two colours.
The faces coloured alike belong to five regular tetrahedra forming the
compound described on page 135, cf.\ Coxeter, Regular polytopes, pp.\ 50,
106.''

Our theorem below asserts that there are 240 colourings of the vertices of the
regular dodecahedron which satisfy our requirement.
The ``paradox'' that this is different from 4 such colourings, as given
in \cite{RBC}, will be resolved in our Corollary. 
We give two proofs of
our theorem: first a combinatorial one, and second a geometrical one.

In our paper we will use the concepts and statements
from \cite{W7}, 
\cite{W1}, \cite{W2}, \cite{W3}, \cite{W4}, \cite{W5}, \cite{W6} ---
the  compound of five regular tetrahedra is defined
in them. We will give its description in {\bf{1.1}}.
Fig.~2 shows the photograph of a cardboard model of the
compound of five regular tetrahedra, taken by A. Lengyel.
(Due to space limitations this figure does not appear in the arxiv version.)




Quite recently
the concept of regular compounds was defined in a strict mathematical
sense by P. McMullen \cite{MM} in 2018, as
follows. A {\it{vertex-regular compound}} consists of copies of a regular
polytope $Q$ with vertex sets in the vertex set of a regular polytope $P$, such
that some subgroup of the symmetry group of $P$ is transitive both on the
vertex set of $P$ and on the copies of $Q$. A {\it{face-regular compound}}
is defined in the dual way, replacing vertices by face-planes. Then a 
{\it{regular compound}} means either a vertex-regular, or a face-regular
compound. There are five regular compounds in ${\Bbb{R}}^3$, cf.\ \cite{MM},
which were
however known and thoroughly investigated already by E. Hess \cite{H} in 1876.
For the compound of five regular tetrahedra, with vertex sets in the vertex set
of the regular dodecahedron, cf. \!\cite{H}, p. \!45.

We also refer to \cite{FT} and \cite{C}
for ample material on regular
figures, and to \cite{W} for its nice figures of cardboard models of regular
polyhedra and of compounds of regular polyhedra. 
At the end of \cite{FT}, German edition, there are
spectacular anaglyphs I-III, displaying these regular
compounds three-dimensionally. 
The compound of five regular tetrahedra is in \cite{FT}, p. 78  --- 
in German p. 79 --- and
in \cite{W}, p. 44 in the 1989 reprint of the first paperback edition in 1974.


\subsection{The compounds of five tetrahedra, and the groups $I$, $I_h$ and
$S_5 \times \{ 1, -1 \} $}

Since we will use only the regular dodecahedron and tetrahedron, and also only
regular compounds,
usually we will drop the attribute ``regular''. We will 
write ``face, edge, vertex'' for those of the dodecahedron.
We will suppose that 
\begin{equation}\label{1}
\begin{cases}
{\text{the centre of the dodecahedron is the origin }}0, \\
{\text{the dodecahedron is inscribed in the unit sphere}} \\
S^2 {\text{ about }}0, {\text{ and has one vertex at the north pole.}}
\end{cases}
\end{equation}

For a graph whose vertices are coloured by some colours, a {\it{colour class}}
is the set of all vertices coloured by a particular colour.

We write $S_5 = S_5(1, \ldots , 5)$ for the {\it{symmetric group}} on 
$\{ 1, \ldots , 5 \}$, 
and $A_5 = A_5(1, \ldots , 5)$ for the {\it{alternating group}} on 
$\{ 1, \ldots , 5 \} $. Further, ${\text{id}}$ is the identical transformation
of ${\Bbb{R}}^3$, and then $-{\text{id}}$ is central symmetry with respect to
the origin. Analogously, we will use the notations $S_5(T_1, \ldots , T_5)$
and $A_5(T_1, \ldots , T_5)$ for the symmetric and alternating group on five
tetrahedra $T_1, \ldots , T_5$. The group $\{ 1, -1 \} $ is written
multiplicatively.

\medskip

\noindent{\bf Definition.}
On the set of all colourings of the vertices of the regular dodecahedron,
such that each face has vertices of all five colours, we introduce the 
action of the group 
$G := S_5 \times \{ 1, -1 \} $ --- i.e., a homomorphism of $G$ to
the permutations of all such colourings --- as follows. It suffices to give
this action for $S_5$ and for $\{ 1, -1 \} $ separately. A permutation $p$
of $\{ 1, \ldots , 5 \} $ acts on a colouring as follows: 
each colour $i$ is
changed to the colour $p(i)$. The action of $1$ is identity, while the action
of $-1$ exchanges the colours of each two antipodal vertices.

\medskip

Observe that here we have a homomorphism from the product group 
$S_5 \times \{ 1, -1 \}$, since the actions of these two groups commute.

We have $I_h = I \times \{ {\text{id}}, -{\text{id}} \} $,
cf.\ \cite{W2} and \cite{W5}  --- their definitions cf. at the beginning of
{\bf{1}}.

By \cite{W3} there are ten tetrahedra with vertex sets in the vertex set of the
dodecahedron. They can be found in Figures 3 and 4: 
\begin{equation}\label{2}
\begin{cases}
{\text{the vertices with the same colours are the vertices of}}\\
{\text{five tetrahedra in Fig.~3, and of other five tetrahedra}} \\
{\text{in Fig.~4. These five tetrahedra, either in Fig.~3 or in}} \\ 
{\text{Fig.~4, together form a compound of five tetrahedra.}} 
\end{cases}
\end{equation}
The two
compounds are congruent, but only via an orientation reversing congruence.
Thus these compounds have a {\it{chirality}}, i.e., left/right-handedness.

Then $I$ acts on these altogether 
ten tetrahedra by permutations. More exactly, 
$I$ acts separately on the five tetrahedra forming any of the two
compounds of five tetrahedra.
Moreover, it acts on any of these two quintuples of tetrahedra
via even permutations, i.e., this action is actually 
a homomorphism $I \to A_5(T_1, \ldots , T_5)$, where $T_1, \ldots , T_5$ are
the tetrahedra in one of our compounds. Moreover, here we have an isomorphism,
so 
\begin{equation}\label{3}
I \cong A_5(T_1, \ldots , T_5) \cong A_5 {\text{ and }}
I_h = I \times \{ {\text{id}}, -{\text{id}} \} \cong A_5 \times \{ 1, -1 \} .
\end{equation}


\section{The theorem and its two proofs}




\vbox{\centerline{\epsfxsize=70mm \epsfbox{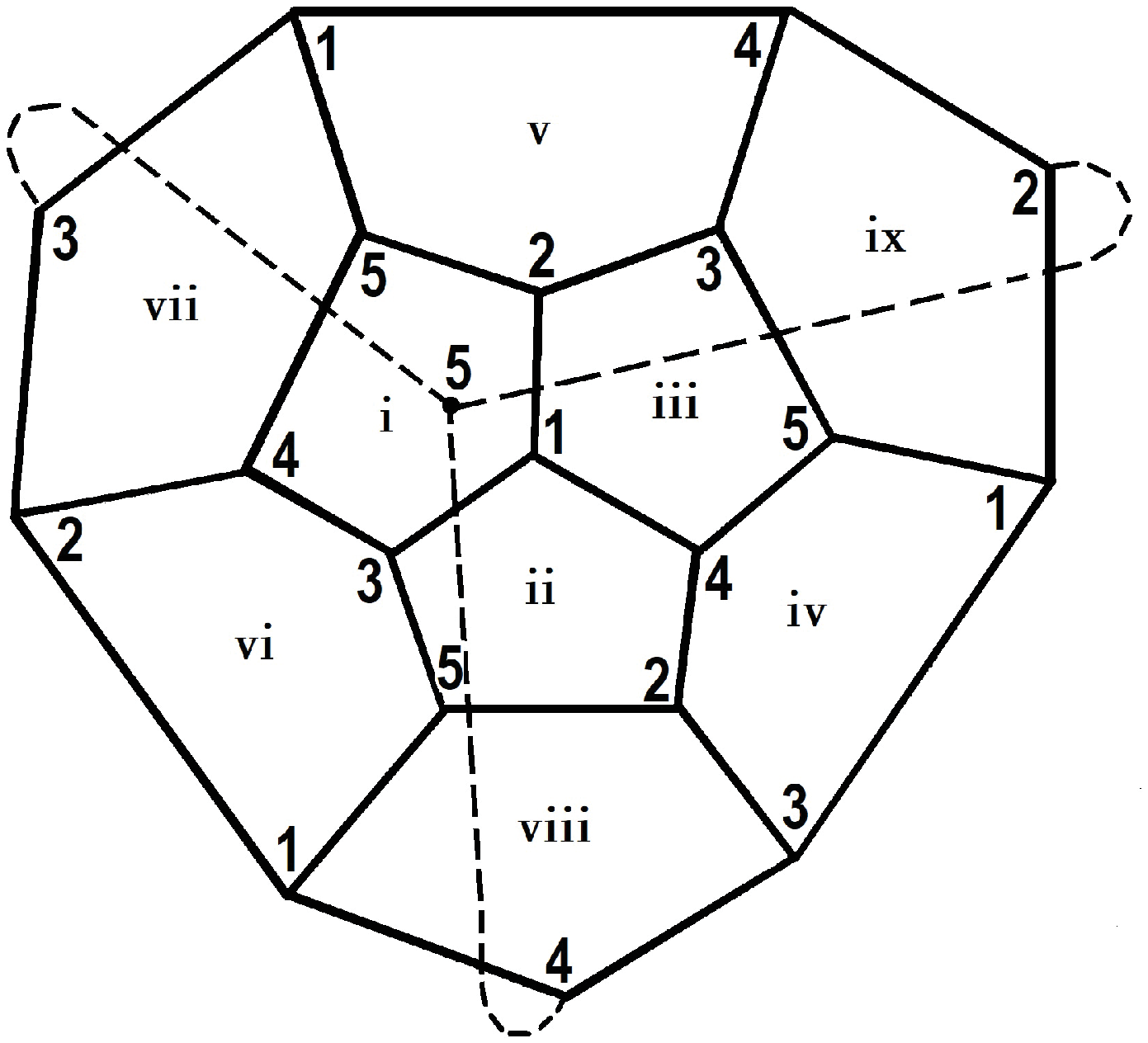}} 

\smallskip
\centerline{
Figure 3.
The first colouring.}
\par}




\vbox{\centerline{\epsfxsize=70mm \epsfbox{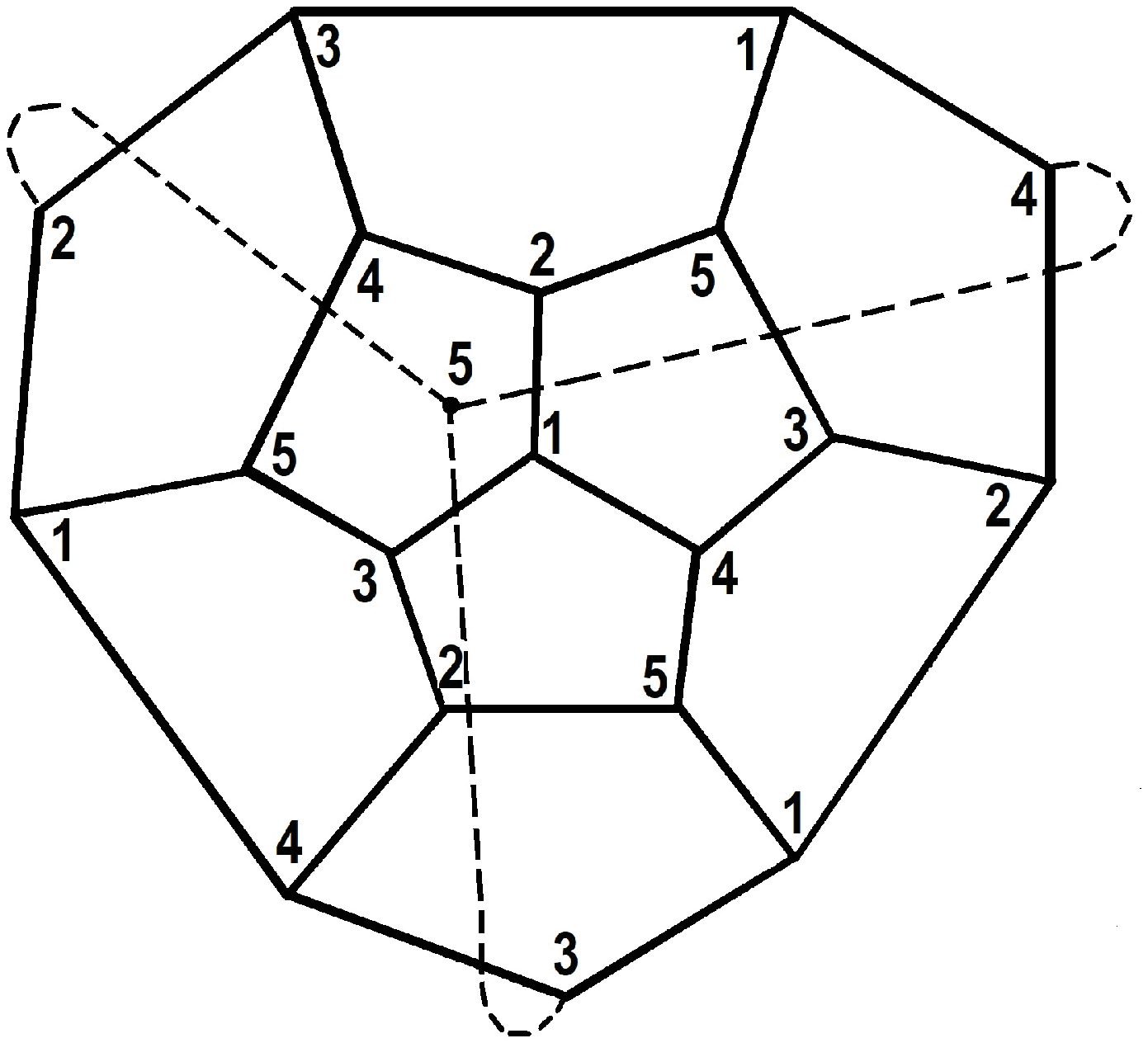}} 

\smallskip
\centerline{
Figure 4.
The second colouring.}
\par}

\smallskip


First we give two colourings of the dodecahedral graph, satisfying our
requirements, in Figures 3 and 4. Recalling (\ref{1}),
in these figures the middle vertex, both times of colour~1, is the north pole
of the unit sphere.
The three neighbours of the north pole, both times of colours 2, 3, 4  --- in 
the same order --- lie on a circle of latitude $C_1$ on the northern hemisphere.
In Fig.~3 the visible faces are denoted by $F_1, \ldots , F_9$.

The numbers of the 2nd, 3rd and 4th neighbours of the north pole are 6, 6
and 3, and they lie on circles of latitude $C_2$, $C_3$ and $C_4$,
respectively.
The unique 5th neighbour of the north pole is the south pole.
$C_1$ and $C_2$ lie on the open northern hemisphere, and $C_3 = -C_2$ and $C_4
= -C_1$ on the open southern hemisphere.
The Euclidean distance between the points of 
$C_3$ and the north pole is the edge-length of the
inscribed regular tetrahedron of the unit sphere.

These colourings have two remarkable properties.

{\bf P1)} For Fig.~3 from any vertex, of colour $i$, we can reach the other
vertices of colour $i$ as follows.
We pass on an edge incident to this vertex, then at
the other endpoint of this edge we turn to the left-hand side edge, and then
at the other endpoint of this second edge we turn to the right-hand side edge:
a ``left-right zigzag''.
For the north pole we pass thus, e.g., through vertices of colours 1, 2, 5, 1.
For any $i \in \{ 1, \ldots ,5 \} $ 
there are exactly four vertices of colour $i$, and they are
the vertices of a tetrahedron, with vertices in the vertex set of
the dodecahedron.
For $i \in \{ 1, \ldots , 5 \}$ 
these five tetrahedra form a compound of five tetrahedra.
The colour of the vertex antipodal to any vertex $v$ is the colour different
from the colours of $v$ and its three neighbours.
For Fig.~4 all above statements of {\bf{P1)}}
hold with exchanging ``left'' and ``right''.
The respective compound of five tetrahedra is obtained from the one in
Fig.~3 by taking its centrally symmetric image with respect to the origin.

{\bf P2)} For Fig.~3, 
on the 12 faces of the dodecahedron, the cyclic orders
of the colours of the vertices, taken in the positive sense,
are exactly all 12 odd cyclic orders. A cyclic order $ijklm$ corresponds
to a permutation $1 \mapsto i, \ldots , 5 \mapsto m$, whose parity is the
parity of the cyclic order. This is well defined, since it is the same as
for the cyclic order
$jklmi$. On opposite faces we have inverse cyclic orders, cf. e.g., $F_4$
and $F_7$.
For Fig.~4 all above statements of {\bf{P2)}} 
hold with changing ``odd'' to ``even''.


\medskip

\noindent{\bf Theorem.}
{\em The number of the five-colourings of the 
vertices of the dodecahedral graph,
such that on each face of the dodecahedron the vertices are coloured with all
five colours, is $240$.
All such colourings can be obtained from any given such colouring, by
the action of some element of the group $G$ from the Definition.
For any such colouring, the colour classes form the vertex sets of 
five regular tetrahedra. These tetrahedra together form any of the
two compounds of five tetrahedra, with vertex sets in the vertex set of
the dodecahedron.} 

\medskip


{\bf First proof of the Theorem.}
Recalling (\ref{1}), we choose 
a permutation of the five colours, such that the vertex of the
dodecahedron at the north pole has colour 1, and its three
neighbours have colours 2, 3, 4, in the order as in Figures 3 and 4.
We will extend the colouring to the vertices on
$C_2$, $C_3$, $C_4$, and the south pole, consecutively.

Both for Fig.~3 and Fig.~4, face $F_1$ 
with one vertex at the north pole, and neighbouring vertices of
colours 2, 3, 
has two yet uncoloured vertices, which must have colours 4, 5.
These can be coloured in two ways: as in Fig.~3, and as in Fig.~4.

We begin with the case shown in Fig.~3. We investigate the six vertices on
$C_2$. Consider the remaining two faces $F_2$, $F_3$, having the north pole as
a vertex, and having neighbouring vertices of colours 3, 4, and 4, 2,
respectively. Both of them have two
yet uncoloured vertices, which must have colours 2, 5, and 3, 5, respectively.
The vertex of $F_1$ of colour 5 enforces that $F_3$ has vertices of cyclic order
$21453$ in the positive sense.
Then the vertex of $F_3$ of colour 5 enforces that $F_2$ has vertices of 
cyclic order $41352$ in the positive sense.
Thus the vertices on $C_2$ are uniquely coloured.

We turn to the six vertices on $C_3$. Each of the faces
$F_4$, $F_5$, $F_6$ has two yet uncoloured vertices, which must have
colours 1, 3, and 1, 4, and 1, 2, respectively. Among
these six vertices on $C_3$, taken in cyclic order, there cannot be two
neighbourly vertices of colour~1.
Therefore each second of them in this cyclic order has colour~1.
By the vertex of colour $4$ of $F_1$, the cyclic order of the vertices of $F_5$
is $52341$ in the positive sense, and then the cyclic orders of the vertices of
$F_6$ and $F_4$ are $53421$, and $54231$, in the positive sense, respectively.
Thus the vertices on $C_3$ are uniquely coloured.

We turn to the three vertices on $C_4$.
Each of the faces $F_7$, $F_8$, $F_9$ 
has a unique yet uncoloured vertex, which therefore must have colours 3, 4,
2, respectively.
Thus the vertices on $C_4$ are uniquely coloured.

At last the south pole has neighbours of colours $3$, $4$, $2$, and second 
neighbours of colours $1$, $2$, and $1$, $3$, and $1$, $4$, respectively.
Therefore its colour must be~$5$.

Thus for Fig.~3 the colouring is uniquely determined, up to a permutation of
the colours.

For Fig.~4 an analogous simple argument shows the same.

Since at the beginning of the proof we have fixed a permutation of the five
colours, both Fig.~3 and Fig.~4 represent in fact $5! = 120$ solutions.

Recalling the Definition, let us
apply an even/odd permutation to a colouring. Then the $120$ colourings 
represented by Fig.~3 or Fig.~4 are taken to colourings
represented by the same figure. Moreover, the parities 
of the cyclic orders of the colours of the vertices on the faces get
preserved/reversed. 
Once more recalling the Definition, let us 
apply $-{\text{id}}$ to a colouring. Then the $120$ colourings 
represented by Fig.~3 or Fig.~4 are
taken to colourings represented by the other figure. Moreover,
``left-right zigzags'' are 
exchanged by ``right-left zigzags''. Observing
that the sets of colour classes in Figures 3 and 4 are different,
the 120 colourings represented by Fig.~3 are different from the 120
colourings represented by Fig.~4. Therefore the total number of colourings
of the dodecahedral graph, satisfying our requirements, is $240$.
These prove the first and second statements of the theorem.

The third statement of the Theorem, i.e., the one about the colour classes,
follows from (2).
\hfill $\square$


\medskip

Observe that both the statements and the proof of our Theorem
are combinatorial. However, one can prove our Theorem also geometrically, as
follows.


\medskip

{\bf Second proof of the Theorem.}
Recalling ({\ref{1}}), any colour class of
a colouring satisfying the hypotheses of our Theorem is a subset of
$S^2$. Moreover,
among its points there occur no smallest and second smallest Euclidean
distances among the vertices of the dodecahedron  --- 
these are the Euclidean
distances between the north pole and the points of $C_1$ and $C_2$. In other
words, {\it{the Euclidean distances between points of any colour class are at
least}} the third smallest Euclidean distance among the vertices of the
dodecahedron. This distance is the Euclidean
distance between the north pole and the points of $C_3$, i.e.,\ {\it{the edge
length of a regular tetrahedron inscribed in $S^2$}}.

However, a subset of
$S^2$ with this last
property has at most four points, with equality only if these
four points are the vertices of a regular tetrahedron, cf.\ \cite{FT}, 
p. 227  --- in the German edition p. 215. Since all
five colour classes have altogether $20$ points, each colour class has exactly
four
points, which are the vertices of a tetrahedron, and also belong to
the vertex set of the dodecahedron. The number of all such tetrahedra is $10$,
and they are depicted in Figures 3 and 4, as tetrahedra with
vertices of the same colours. 
Then all five colour classes are
either as depicted in Fig.~3, or as depicted in Fig.~4, or some colour class
is as depicted in Fig.~3 and some other colour class is as depicted in
Fig.~4. Since the colour classes are disjoint, from Figures 3 and 4
we see that in the last case any two such
colour classes are centrally symmetric images of each other with respect to
$0$. Hence the number of colour classes is two, a contradiction. Therefore the
colour classes 
are either as depicted in Fig.~3, or as depicted in Fig.~4. Now one can finish 
this proof as the combinatorial proof. 
\hfill $\square $

\medskip


By our Definition,
any element of $G$ takes a colouring satisfying our requirements to
another such colouring. Also conversely, by our Theorem,
any such colouring can be taken over
to any other such colouring by applying some element of $G$. Therefore
we can say that, up to the action of elements of $G$,
we have ``geometrically'' a unique colouring.

If $H$ is an arbitrary subgroup of $G$, we can ask for the {\it{number of our
colourings, up to the action of an element of $H$}}.
This number is $|G|/|H| = 240/|H|$.
For the statement of W. W. Rouse Ball -- H. S. M. 
Coxeter \cite{RBC}, p. 242, we consider 
$H := I \times \{ {\text{id}} \} \cong A_5 (T_1, \ldots , T_5)
\times \{ {\text{id}} \} \cong A_5 \times \{ 1 \} $, cf.\ ({\ref{3}).
Then we have the following


\medskip

\noindent{\bf Corollary.} (\cite{RBC}, p. 242)
{\em The number of the colourings of the dodecahedral graph, satisfying the
requirements of our Theorem, 
up to the application of any element of $ I \times \{ {\text{\rm{id}}} \}
\cong A_5 \times \{ 1 \}$, 
is $4$. They can be distinguished as follows.

1) The colour classes are the sets of vertices of five regular
tetrahedra, forming
either of the two compounds of five regular tetrahedra, with vertex sets in
the vertex set of the regular
dodecahedron. These two possibilities are depicted in Figures 3 and 4, and
correspond to ``left-right zigzags'' or ``right-left zigzags''.

2) Further they can be distinguished by odd/even
cyclic orders of the colours of the vertices on
one (or all) face(s) of the regular dodecahedron
in Figures 3 and 4.

These two choices being independent, we have altogether four possibilities.} 
\hfill $\square$


\medskip

{\bf{Acknowledgement.}} We thank A. Lengyel for his photograph of the compound
of five tetrahedra, and for his kind courtesy to agree to
include it in our paper. We also express our gratitude to the anonymous
referee, whose suggestions have greatly improved the presentation of our
material.


\enddocument